\documentclass{birkart}
\usepackage{graphics}

 \newtheorem{thm}{Theorem}[section]

 \newtheorem{cor}[thm]{Corollary}

 \newtheorem{lem}[thm]{Lemma}

 \newtheorem{prop}[thm]{Proposition}

 \theoremstyle{definition}

 \newtheorem{defn}[thm]{Definition}

 \theoremstyle{remark}

 \newtheorem*{ex}{Example}

 \numberwithin{equation}{section}

\DeclareSymbolFont{AMSb}{U}{msb}{m}{n}
\DeclareMathSymbol{\N}{\mathbin}{AMSb}{"4E}
\DeclareMathSymbol{\Z}{\mathbin}{AMSb}{"5A}
\DeclareMathSymbol{\R}{\mathbin}{AMSb}{"52}
\DeclareMathSymbol{\Q}{\mathbin}{AMSb}{"51}
\DeclareMathSymbol{\I}{\mathbin}{AMSb}{"49}
\DeclareMathSymbol{\C}{\mathbin}{AMSb}{"43}

\begin{document}

\title[Soap films]
{Cartan's magic formula and soap film structures}

\author[J. Harrison]{Jenny Harrison}

\address{  Department of Mathematics\\ University of California, Berkeley\\Berkeley, CA\\
94705\\USA}

\email{harrison@math.berkeley.edu}



\begin{abstract}
 
A soap film is actually a thin solid fluid bounded by two surfaces of opposite orientation.   It is natural 
to model the film using one
polyhedron for each side. Two problems are to get the polyhedra for both sides
to be in the same place without canceling each other
out and to model triple junctions without introducing extra boundary components.
We use chainlet geometry to create dipole cells and mass cells which accomplish these goals and model
faithfully all observable soap films and bubbles.     
We introduce a new norm on chains of these cells and prove lower semicontinuity of area.  A geometric version of
Cartan's magic formula provides the necessary boundary coherence.

\end{abstract}

\maketitle

\section{Introduction}
Plateau observed that soap films have only two possible kinds of branching:  (1) three sheets of surface meeting at
$120^{\circ}$ angles along a curve and (2) four such curves meeting at approximately  $109^{\circ}$ angles at a point.  
The question known as {\em Plateau's Problem} naturally arose,
\begin{quotation} Given a loop of wire in 3-space, is there a surface with minimal area spanning it?
\end{quotation} 

 Osserman \cite{oss} wrote  ``the question of how to formulate this problem precisely has been
almost as much of a challenge as devising methods of its solution.''
 Over the past century several mathematical models  for area minimizing surfaces have been produced. The
regularity and singularity structure of the resulting surfaces have been extensively studied,  but none are sufficiently
general to include all types of surfaces that arise as soap films.

The first solution to Plateau's Problem was due to J. Douglas \cite{Douglas} who established the
existence of area-minimizing disks for a given Jordan curve.    He  showed
that amongst all mappings
of the disk into $\R^3$ such that the disk boundary is mapped
homeomorphically onto a given
curve $\gamma$, there exists at least one mapping whose image has smallest
possible area.    Douglas was awarded one of the first Fields' 
medals for this work.   

Federer and Fleming \cite{FF} proved the existence of area-minimizing embedded orientable
surfaces.  They found an integral current with minimal mass spanning a given
curve.  Fleming later proved  that this integral current was an embedded,
orientable surface.     
Almgren used varifolds \cite{A2} to treat the nonorientable case but
found no natural boundary operator.   Soap film regularity of varifolds remains an open question.
Brakke \cite{Brakke} produced interesting models of soap films using a
covering space of the complement of the soap film boundary.   There remain questions about existence and regularity of
area minimizers in this category and it is not known whether all soap films are modeled.

Using methods and results from chainlet geometry (\cite{continuity}, \cite{iso})  summarized in $\S 2,$ we define new models
for soap films called {\em dipolyhedra} in $\S 3$, giving new definitions to the notions of {\em curve, surface,} and  {\em area}.   In
$\S 4$ we prove lower semicontinuity of area and in $\S 5$ we show that 
dipolyhedral models apply to all soap films observed by Plateau,  soap bubbles\footnote{The author would like to thank Harrison Pugh
for the artwork in Figures 1 and 3.}, the surfaces studied by Douglas, Federer and Fleming, as well as  films with singular branched
curves,  nonorientable films, films touching only part of the Jordan curve boundary and other examples that have eluded previous
models.

\begin{figure}\label{fig.1}
 \vspace{2.5in}
 \hspace{-3.0in}
%\special{picture bubbles.ps}
\begin{center} 
\resizebox{4.0in}{!}{\includegraphics*{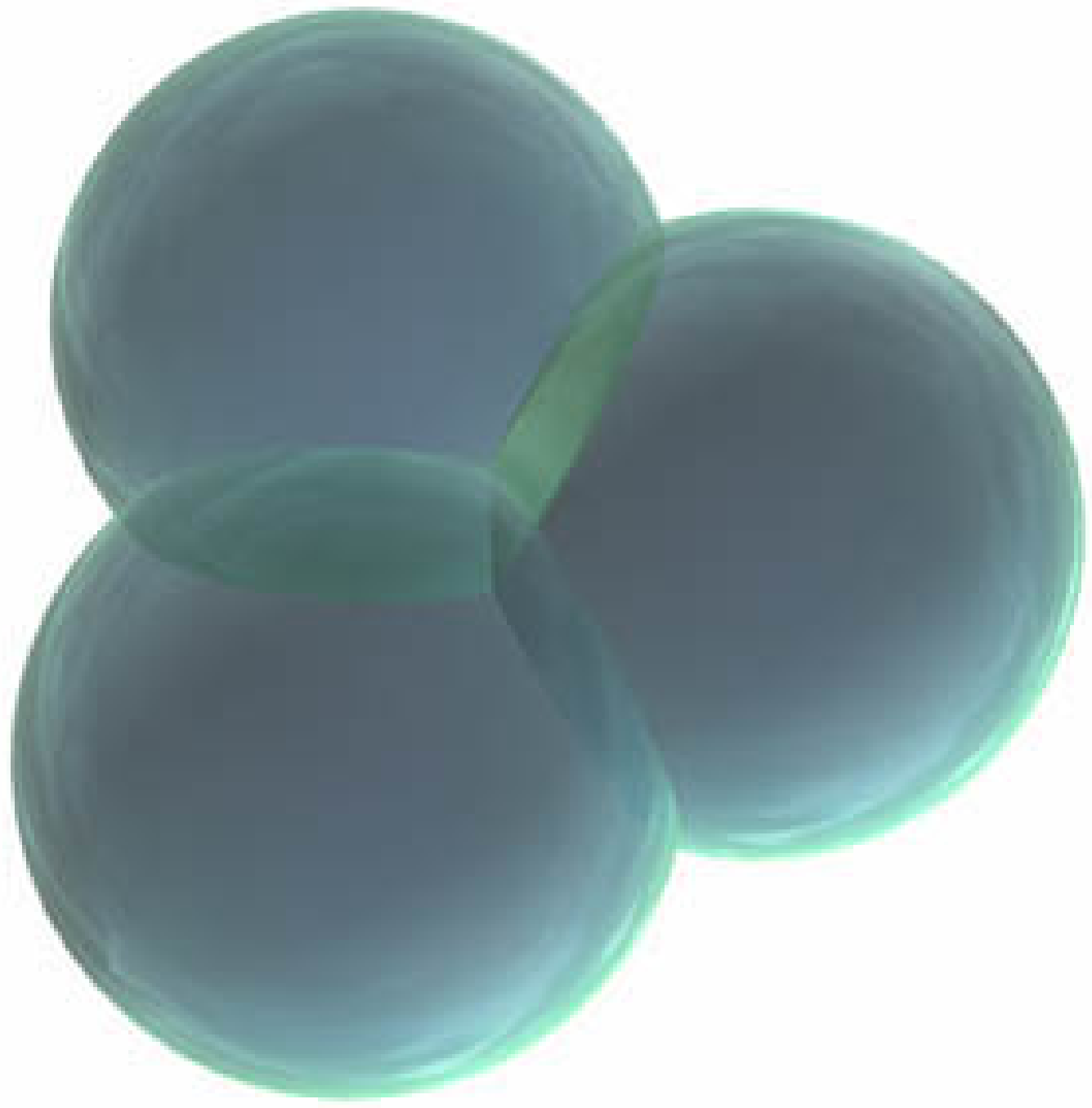}} 
\end{center}
\caption{\quad}
\end{figure}

 \section{Normed groups of polyhedra}
 
The concept of continuity of domains plays an important role in our approach.  In the first
two sections we  develop the concept of  {\em closeness} for domains.   This is conveniently done by starting with 
norms on polyhedra.   

For
$k
\ge 1$, a
$k$-dimensional {\em   cell} in $\R^k$ is defined to be the finite intersection of open half spaces.
A $k$-cell   in $\R^n$ is a subset of a $k$-plane $\Pi$
in
$\R^n$ that is a
$k$-cell in $\Pi.$   All   $k$-cells are assumed to be oriented.
A $0$-cell is a single point $\{x\}$ in $\R^n$.  No
orientation need be assigned to the
 $0$-cell.    One may take formal sums of   $k$-cells with integer coefficients $S =
\sum a_i
\sigma_i
$ and form equivalence classes of
 polyhedral chains $P = [S]$.   
Two formal sums are equivalent if integrals of smooth differential forms are the same over each of
them.    Let
$\mathbf{P}_k(\Z)$ denote the abelian group of 
$k$-dimensional polyhedral
$k$-chains in $\R^n$ with integer coefficients.    Let $\partial$ denote
the usual boundary operator on
 cells.   For $k = 0$ we set
$\partial P = 0.$  The  sequence of polyhedral chains of arbitrary dimension and
$\partial$ form a
chain complex
$\mathbf{P}_*(\Z).$

  For $k \ge 1$, let
$M(\sigma)$ denote the 
$k$-dimensional Hausdorff measure  of a $k$-dimensional  cell $\sigma$.
Every $0$-dimensional
   cell  $\sigma^0$ takes the form
$\sigma^0 =
\{x\} $ and we  set
$M(\sigma^0) = 1.$
The {\em mass} of $P= \left[\sum a_i \sigma_i\right]$ is defined by

$$M(P) = \sum_i|a_i|M(\sigma_i)$$  where the $\sigma_i$ are
non-overlapping.  Mass is a norm on the   group  $\mathbf{P}_k(\Z).$

 \begin{defn}[Multicellular chains]    If
$\sigma$ is a cell and $v$ a vector let $T_{v}\sigma$ denote the translation of
$\sigma$ through  
$v$.

     A {\em
$1$-multicell} is a cellular chain of the form

$$\sigma^1 = \sigma^0-T_{v_1} \sigma^0 $$  where $\sigma^0$ is a cell and $v_1$ is a vector.
We inductively define
  $2^j$-multicells. Given a    vector    $v_j$ and a
$2^{j-1}$-multicell $\sigma^{j-1}$, define the {\em   $2^j$-multicell}
$\sigma^j$ as the cellular chain $$\sigma^j = \sigma^{j-1} -T_{v_j}\sigma^{j-1}.$$  Thus
$\sigma^j$ is generated by  vectors $v_1, \dots ,v_j$ and a cell $\sigma^0.$
 
\medskip An {\em  integral $2^j$-multicellular  chain in $\R^n$} is a formal sum of
$2^j$-multicells, $S^j =
\sum_{i=1}^n a_i \sigma_i^j $ with coefficients $a_i \in \Z$.   Let $[S^j]$ denote the polyhedral chain represented
by $S^j.$  

  Given a $2^j$-multicell $\sigma^j$ generated by a cell $\sigma^0$ and vectors   $v_1, \cdots, 
v_j$,  define $\|\sigma^0\|_0 = M(\sigma^0) $, $ \|\sigma^1\|_1 =  M(\sigma^0)\|v_1\| $ and for $j \ge 2$,
$$\|\sigma^j\|_j = M(\sigma^0)\|v_1\|\|v_2\|\cdots  \|v_j\|.$$ 
 For $S^j = \sum a_i \sigma_i^j$ define $$\|S^j\|_{j} = \sum_{i=1}^n
|a_i|\|\sigma_i^j\|_{j}.$$

 \end{defn}

\begin{defn}[Natural norms]
Let $P \in {\mathbf P}_k(\Z)$ be a polyhedral $k$-chain  and $r \in
\Z^+$.  For $r = 0$ define
$$|P|^{\natural_0} = M(P).$$
For $r \ge 1$ define the {\em $r$-natural} norm
\begin{equation}\label{natural}
                 |P|^{\natural_r} = \inf\left\{ \sum_{s=0}^r \|S^j\|_j   +
  |C|^{\natural_{r-1}}\right\}
\end{equation}

where the infimum is taken over all  decompositions
                                                                          $$P
= \sum_{s=0}^{r} [S^j]  + \partial C $$ where $S^j$ is a $2^j$ multicellular $k$-chain and $C$ is a
polyhedral
$(k+1)$-chain.    
\end{defn}

The  $r$-natural norm was
introduced     in \cite{iso} for real coefficients and $\R^n.$   In this paper we denote this norm by 
$|P|_{\R}^{\natural_r}$. 

\begin{lem}  $  |P|^{\natural_r} $ is a norm on the group of polyhedral chains $\mathbf{P}_k(\Z)$.  
\end{lem}

\begin{proof}   It is clear that  $  |P|^{\natural_r} $  is a semi-norm with 
$|P|_{\R}^{\natural_r}
\le |P|^{\natural_r}.$  If
$|P|^{\natural_r} = 0$ then
$|P|_{\R}^{\natural_r} = 0$ and hence $P = 0$ since $|P|_{\R}^{\natural_r}$ is a norm.     
\end{proof}

 The group of polyhedral $k$-chains ${\mathbf P}_k(\Z)$ completed with
the norm $|\quad|^{\natural_r}$ is denoted ${\mathbf
N}_{k}^r(\Z).$ The elements of ${\mathbf N}_{k}^r(\Z)$ are
called {\em $k$-dimensional chainlets of class $N^r$.
}        
  The flat norm of Whitney \cite{W} has played an important role in
geometric measure
theory as flat integral chains correspond to integral currents and thus flat integral chains give a geometric approach to the study of integral currents.  
\cite{FF}
\cite{Fleming2}.
  Chainlets provide a geometric approach to the study of {\em integrable} currents including those with unbounded
mass \cite{currents}.  Unlike currents, they permit coefficients in an abelian group.  In the next section we introduce
two important examples of chainlets,   {\em mass} and {\em dipole} cells.

It follows immediately from the definitions that the boundary
operator on polyhedral chains is
bounded w.r.t. the
$r$-natural norms.

\begin{lem} \label{lem.boundary}
If $P \in {\mathbf P}_k(\Z)$ then

$$|\partial P|^{\natural_{r+1}} \le |P|^{\natural_r}. $$
\end{lem}

Therefore the boundary $\partial A$ of a $k$-dimensional chainlet $A$ of class $N^r$ has unique
definition as a $(k-1)$-dimensional chainlet   of class $N^{r+1}.$ 
  
\begin{defn}[Support of a chainlet]   The support of a cell $\sigma$ is merely the points contained in its closure
$|\sigma|.$  If
$P$ is a polyhedron it can be written $P = \sum a_i \sigma_i$ where the $\sigma_i$ are non-overlapping.  The support of $P$
is defined as 
$|P| = \cup  
|\sigma_i|.$   We say a closed set $F$ {\em supports} a chainlet
 $A$ if for every open set $U$
containing $F$ there is a
sequence $\{P_j\}$ of  polyhedra tending to $A$ such
that each $P_j$ lies in $U$.  If there is a smallest set $F$
which supports $A$ then $F$ is called the {\em support} of $A$ and
denoted $| A|$.
\end{defn}

 \newpage
\section{Normed groups of dipolyhedra}\label{sec.chains}  \setcounter{equation}{0}
 
  We now specialize to two examples of chainlets used to model soap films.  
 
\subsubsection*{Mass chains}

 Suppose
$\tau$ is a 
$(k-1)$-cell in $\R^n$,
$0
\le k-1 \le n$ and
$v$ is a vector in $\R^{n+1}$. Let  
$\tau \times 2^{-j}v$ denote the oriented orbit of  
$\tau$  through the time-$t$ map of  the vector field $v(x) = v$ for $0 \le t \le 2^{-j}.$
 (Here $2^{-j}v$ refers to scalar multiplication
of a vector $v$ by   $2^{-j}$. )         If $v$ is transverse to $\tau$ then $\tau \times 2^{-j}v$
is a $(k+1)$-cell.  Otherwise it is degenerate.

\begin{lem} \label{massnat} The sequence of weighted cells $2^j(\tau \times 2^{-j}v)$
  is Cauchy in the $1$-natural norm.    
\end{lem}

\begin{proof}     Note that
$$\tau \times 2^{-j}v = (\tau \times 2^{-j-1}v) + (T_{ 2^{-j-1}v}\tau \times  2^{-j-1}v).$$
Since 

\[
\begin{array}{rll} \frac{\tau \times 2^{-j}v}{2^{-j}} -\frac{\tau \times 2^{-{j-1}}v}{2^{-{j-1}}} &=  \frac{(T_{
2^{-j-1}v}\tau \times  2^{-j-1}v) -(\tau \times 2^{-j-1}v) }{2^{-j}} \\&= \frac{T_{ 2^{-j-1}v}(\tau \times 2^{-j-1}v)
  -(\tau \times 2^{-j-1}v) }{2^{-j}} 
\end{array}
\]  
we have
\[
\begin{array}{rll}  \left| \frac{\tau \times 2^{-j}v}{2^{-j}} -\frac{\tau \times 2^{-{j-1}}v}{2^{-{j-1}}} 
\right|^{\natural_1}
&=\left|  \frac{T_{ 2^{-j-1}v}(\tau \times 2^{-j-1}v)
  -(\tau \times 2^{-j-1}v) }{2^{-j}}\right|^{\natural_1} \\&\le 
\|v\|^2|\tau|^{\natural_0}2^{-j}/4.
\end{array}
\]  Thus $2^j(\tau \times 2^{-j}v)$ is Cauchy in the $1$-natural norm.
\end{proof}

 We denote the limit in the $1$-natural norm by $$\mu_v \tau = \lim_{j \to \infty}
2^j(\tau \times 2^{-j}v)$$ and call this a  {\em  directional mass cell in the direction $v$}.  If $v$ is  unit
and normal to
$\tau$ then
$\mu_v \tau$ is called a  {\em mass cell}.   The {\em mass} of a mass cell $\mu_v
\tau$ is defined to be $M(\mu_v \tau) =   M(\tau).$ 
Observe that 
  $| \mu_v \tau| = | \tau|,$ 
 
An {\em integral mass chain} $T$  is a formal sum $T   = \sum a_i \mu_{v_i}\tau_i$ where $a_i \in \Z.$  It is {\em
nonoverlapping} if the cells $\{\tau_i\}$ are nonoverlapping.    The group of integral mass chains of dimension $k$ is
denoted $T_k(\Z).$ 
 The boundary of $T$  is naturally defined by $\partial T = \sum a_i
\partial \mu_{v_i} \tau_i .$

\begin{ex}
 A Dirac delta distribution supported at a point $p \in \R$ corresponds to a  $1$-dimensional mass cell 
$\mu_{e_1}\{p\}
$ where
$e_1 = 1
\in
\R^1$  in the sense that if $f$ is a Lipschitz function then $\delta_p(f) = \int_{\mu_{e_1}\{p\}} f = f(p).$
\end{ex}
\begin{ex}
A   Moebius strip supports a $3$-dimensional mass chain.  
 
\end{ex}

\subsubsection*{Dipole chains}  Let $\sigma$ be a $k$-cell.  For $v$ transverse to $\sigma$, define $$\delta_v
\sigma =\mu_v
\partial
\sigma +\partial
\mu_v
\sigma .$$ Since $\mu_v
\sigma$   is a $(k+1)$-dimensional chainlet   of class $N^1$ then $\partial \mu_v \sigma$ is a $k$-dimensional
chainlet   of class $N^2.$   Since $\mu_v \partial
\sigma $ is a $k$-dimensional chainlet   of class $N^1$ it is also a
$k$-dimensional chainlet   of class $N^2$.  Thus
$\delta_v
\sigma$  is a
$k$-dimensional chainlet   of class $N^2$ with  $| \delta_v \sigma| \subset  | \sigma|$. 

It follows   from the definitions    that the geometric {\em Lie derivative of
$\sigma$ in the direction
$v$} is well defined in the $2$-natural norm
  $$\mathcal{L}_v  \sigma = \lim_{j \to \infty}
2^j(\sigma -T_{2^{-j}v}\sigma).$$  This leads immediately to the next result.

\begin{thm}[Geometric Cartan's magic formula]     $$\mathcal{L}_v  \sigma = \mu_v
\partial
\sigma + \partial \mu_v \sigma .$$ 
\end{thm}   (In a sequel, \cite{discrete}, the author shows how to extend this to Lipschitz vector
fields
$v$ and Lipschitz cells $\sigma$.) 

If $v$ is unit and  normal to
$\sigma$ we call $\delta_v\sigma$ 
 a {\em  dipole cell}.

 Since a dipole cell $\delta_v \sigma$ is a  $k$-dimensional chainlet of class $N^2$, its boundary is a uniquely defined $(k-1)$-dimensional chainlet
of class $N^3$  and is a sum of dipole cells.     It follows immediately from the definitions that $$\partial \delta_v 
\sigma = \delta_v  \partial \sigma .$$ 
\begin{ex} A  dipole point of physics based at $p \in \R$ corresponds to a  $0$-dimensional dipole cell
$\delta_{e_1}\{p\}$.   It is the  boundary of the Dirac delta mass cell $\mu_{e_1}\{p\}.$     
\end{ex}

 Define the {\em weight} of a dipole cell
$\delta_v \sigma$ to be
$W(\delta_v \sigma) = M(\sigma).$   (We do not call this {\em mass} since the $0$-natural norm, or mass, of
$\delta_v \sigma$ is infinite.   ) 

An {\em integral dipole chain} $S$  is a formal sum $S   = \sum a_i \delta_{v_i}\sigma_i$ where $a_i \in \Z.$   A 
dipole chain is {\em nonoverlapping} if the cellular chain
$\sum a_i
\sigma_i$ is nonoverlapping.   
The group of integral
dipole chains of dimension $k$ is denoted $S_k(\Z).$ 
 The boundary of $S$  is naturally defined by $\partial S = \sum a_i
\delta_{v_i}\partial \sigma_i .$

\begin{ex} Consider three one-dimensional  dipole cells in the plane meeting in a point
$p$ forming angles of
$120^{\circ}.$  Their sum $D$ is a   dipole chain with $\partial D$  supported in three
points not including
$p$.  Contrast this with three oriented one-dimensional cells  with  support the same as $|D|$.   Under $\Z$ or
$\R$  coefficients, the support of the boundary must contain $p$, no matter how the orientations are chosen.   See Figure
2.
\end{ex} 

\begin{figure}\label{fig.2}
\vspace{1.5in}
\hspace{-4.0in}
%\special{picture dipole.ps}
\begin{center} 
\resizebox{4.0in}{!}{\includegraphics*{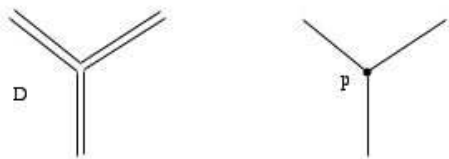}} 
\end{center}
\caption{The support of $\partial D$ in $\Z_2$ coefficients does not include $p$}
\end{figure}

Henceforth, for our applications to soap films, we assume that each mass and dipole cell is supported in $\R^3$ and
the vector $v= e_{4} = (0, 0, 0, 1) \in
 \R^{4}$.   We     simplify the notation and write 
$$\delta \sigma = \delta_{e_{4}}\sigma$$ for a $k$-dimensional dipole cell and
$$\mu \tau = \mu_{e_{4}}\tau$$ for a $k$-dimensional mass cell where $\sigma$ is a   $k$-cell and $\tau$ is a  
$(k-1)$-cell.   If $P = \sum a_i \sigma_i$ is a polyhedron we define $\delta P = \sum a_i \delta \sigma_i$ and $\mu P =
\sum a_i \mu \sigma_i.$   
 
We next permit coefficients in an abelian group
$G$ provided with a translation invariant metric.  Let
$|g|$ denote the distance between $g$ and the identity $0$.  Assume that
$G$ is a complete metric space.    If $H \subset G$ is a
closed subgroup we may put the quotient metric $|\bar{g}| =
\inf\{|g|:g \in  \bar{g}\}$ on
$G/H.$  If $G = \Z$ we take for $|g|$ the usual absolute value.  The group $\Z_2 = \Z/2\Z$ is of special interest.   We
give it the quotient metric.  

Define
$$T_k(G) = G \otimes T_k(\Z).$$  This is the  group of
all mass $k$-chains supported in $\R^n$ with coefficients in $G.$
     If $T \in T_k(G)$ then $T = \sum g_i \mu \tau_i$ where $g_i \in
G $  and the $\{\tau_i\}$    are nonoverlapping 
$(k-1)$-dimensional     cells.  
The mass of  
$T \in T_k(G) $ is defined by
$$M(T) = \sum |g_i| M(\tau_i)$$ where $T =   \sum g_i\mu_{v_i}\tau_i$ is nonoverlapping.  It is easy to verify that
$M$ is a norm on the group of  mass chains $T_k(G).$  
 
Similarly, define   $$S_k(G) = G \otimes S_k(\Z).$$  This is the group of
 dipole chains with coefficients in $G$.   If $S \in S_k(G)$ then $S = \sum g_i \delta \sigma_i$ where $g_i \in
G $  and the $\{\sigma_i\}$    are nonoverlapping 
$k$-dimensional     cells.  
 Let
$S \in S_k(G) $ be a dipole chain with   $S = \sum  g_i\delta_v \sigma_i$  non-overlapping.   The {\em weight} of $S$ is defined by
$$W(S) = \sum |g_i| M(\sigma_i).$$    Then
$W$ is a well defined norm on the group of dipole chains.  This quantity corresponds to the  $k$-dimensional Hausdorff
measure of
$| S|$ with multiplicity.

 \subsection*{Dipolyhedra}
 
  Let $\mathcal{H}^k$ denote $k$-dimensional Hausdorff measure.  
  Define
 $$D_k(G) = S_k(G) \oplus T_k(G).$$
 Define $W(D) = W(S)$ where $D = S+T.$  This is clearly a semi-norm on the group of dipolyhedra.
\begin{lem}  Let $D \in D_k(G)$.  Then $$\mathcal{H}^k(|D|) = W(D).$$
 
\end{lem}   

\begin{lem}\label{Dip}  If $D = S+T$ satisfies $\partial D = 0$ then $|T| = |\partial S|$. 
\end{lem}

\begin{proof} $\partial S = -\partial T$ implies $|\partial S| = |\partial T| = |T|.$
\end{proof}

  Dipolyhedra are  the basic building blocks for soap films  
and fundamental to what follows.    They  carry the simplicity of $\Z_2$ coefficients, they may have triple branch
junctions, yet their boundaries can be supported in a given wire.    The primary role of mass cells is to correct boundary
errors which would arise if one were to model soap films using mod 2 polyhedra.

The {\em support} of a nonoverlapping dipolyhedron $D = \sum  \delta \sigma_i + \sum \mu
\tau_i$ is the same as the support $|P|$ of
the   polyhedron $P =
\sum
\sigma_i +
\sum
\tau_i.$  In this paper, all dipolyhedra satisfy $|T| \subset |S|$ so that $| D| =|S| .$

\begin{defn}

Define the {\em energy norm} of a dipolyhedron
$D = S+T $ to be
$$E(D) = W(S) + M(T).$$ 

\end{defn}

\subsection*{The part of a dipolyhedron in a half space}

Let $H_s$ denote the open half space $$H_s = \{x: x^i < s \mbox{ for
some } 1 \le i \le 3\}.$$  If $\sigma$ is a $k$-cell
  then  $\sigma \cap H_s$ is also a
$k$-cell.     
  Given a dipole cell $\delta \sigma$, define
the  dipole cell
$$(\delta \sigma)
\cap H_s =
\delta (\sigma \cap H_s).$$ Given a mass cell $\mu \sigma$, define
the  mass cell
$$(\mu 
\sigma)
\cap H_s =
\mu (\sigma \cap H_s).$$

Extend the definition to $S \cap H_s$ and $T \cap H_s$ for   dipole chains $S$ and mass chains $T$ by linearity.   If $D =
S+T$ define
$D
\cap H_s = S
\cap H_s + T
\cap H_s$ where $S \cap H_s $ is a dipole chain and $T
\cap H_s$ is a mass chain.\footnote{  In \cite{Fleming2} Fleming used the notation $A \cap X$, Federer used $A
\lfloor{X}.$  We prefer Fleming's notation since it coincides with
the notation of intersection of polyhedral chains with
polytopes and reserve the notation $T \lfloor X$ for currents and forms. }  

 \begin{lem} \label{Qlemma} If $D = S+T$ is a dipolyhedron and $H_s$ is a half space   
$$E(D\cap H_s)
\le E(D).$$
 
If $Q$ is a $3$-cube  then   $$E(D\cap Q)
\le E(D).$$
 
\end{lem}

\begin{proof} Note 
\[
\begin{array}{rll}
W(\delta \sigma \cap H_s)   = 
   M(\sigma \cap H_s)  \le  M(\sigma).   
\end{array}
\]   Suppose $S = \sum  \delta \sigma_i$ is nonoverlapping.     Then   
$$W(S\cap H_s) \le \sum M(\sigma_i) = W(S). $$ 
Similarly,   \[
\begin{array}{rll}
M(\mu \sigma \cap H_s)   = 
   M(\sigma \cap H_s)  \le  M(\sigma).   
\end{array}
\]   Suppose $T = \sum  \mu  \sigma_i$ is nonoverlapping.     Then   
$$M(T\cap H_s) \le \sum M(\sigma_i) = M(T). $$ 
It follows that $$E(D\cap H_s) = W(S\cap H_s) + M(T\cap H_s) \le  W(S) + M(T) = E(D).$$
Observe that $D\cap Q$ is
well defined for a dipolyhedron  $D$ and an  open $3$-cell  $Q $ since $D \cap H_s$ is well defined for every half space $H_s$ and
$Q$ is a finite intersection of half spaces.  Hence the estimates follow for cells $Q$ by writing $Q$ as a finite
intersection of half spaces and applying the previous estimates repeatedly.
\end{proof} 

 For a cell $\sigma$, define the {\em slice}
$\sigma_s = \partial (\sigma \cap H_s) - (\partial \sigma) \cap
H_s.$  Let $u$ be normal to $\sigma$.  Since $\sigma_s$  is a cellular chain 
$\delta (\sigma_s)$ and $\mu  (\sigma_s)$ are a uniquely defined dipole chain and mass chain, resp.

Define  the {\em slice} of a dipolyhedron $D$ by
$$D_s = \partial (D \cap H_s) - (\partial D) \cap H_s.$$   Observe that $D_s$ is a dipolyhedron since $D$ is a dipolyhedron.
 It follows directly from the definitions and relations that $\delta (\sigma_s) = (\delta \sigma)_s $ and $\mu  (\sigma_s) = (\mu \sigma)_s. $
Thus we may write unambiguously $\delta \sigma_s$ and $\mu_s \tau_s$ for the slice of a dipole cell and a mass cell.
 
 \begin{lem}[Fubini inequalities for   dipolyhedra] \label{lintegral}Let  $D$ be a dipolyhedron.  Then  
 $$\int_{-\infty}^{\infty} W(S_s) ds \le W(S),$$   
 $$\int_{-\infty}^{\infty} M(T_s) ds \le  M(T)$$  and
$$\int_{-\infty}^{\infty} E(D_s) ds \le  E(D).$$      
\end{lem}
 
\begin{proof}  
First note the elementary estimate  $$\int_{-\infty}^{\infty} M(\sigma_s)ds \le M(\sigma)$$  
Then
\[
\begin{array}{rll} \int_{-\infty}^{\infty} W(\delta \sigma_s) ds  = \int_{-\infty}^{\infty} M(\sigma_s)ds \le   M(\sigma)
= W(\delta \sigma).
\end{array}
\]

 Suppose $S = \sum  \delta \sigma_i$ is a non-overlapping dipole chain.   Then
$S_s = \sum  \delta \sigma_{is} $  is a non-overlapping dipole chain and 
\[
\begin{array}{rll}
\int_{-\infty}^{\infty}W(S_s)ds   =  \sum  \int_{-\infty}^{\infty}M(\sigma_{is}) 
 \le \sum M(\sigma_i)     =  W(S)  .
\end{array}
\]

The inequality $\int_{-\infty}^{\infty} M(T_s) ds \le
 M(T) $ is similar.   Finally $$\int_{-\infty}^{\infty} E(D_s) ds=  \int_{-\infty}^{\infty} W(S_s) ds  +
\int_{-\infty}^{\infty} M(T_s) ds
\le W(S) + M(T) = E(D).$$ 
  
\end{proof}

\begin{defn} 

For $D \in D_k(G)$ define
$$E_{\flat}(D) =  \inf\{E(B) + E(C)  : D = B + \partial C, B \in D_k(G), C \in 
D_{k+1}(G)\}.$$   
\end{defn}

It is clear that $E_{\flat}$ is a seminorm.  The fact that it is a norm will be established soon.
We complete the group ${D}_k(G)$ with $E_{\flat}$ to obtain
a group $\mathbf{D}_k(G)$ of {\em flat dipolyhedra}.

 \begin{lem}   If $D$ is a dipolyhedron $E_{\flat}(\partial D) \le
E_{\flat}(D)$. 
\end{lem}

\begin{proof}   
Given $\epsilon > 0$, choose $B, C$ so that $D = B + \partial C$ and  $E(B) + E(C) < E_{\flat}(D) + \epsilon.$ Since $\partial D = \partial B$ we have $$E_{\flat}(\partial D)
\le E(B) < E_{\flat}(D) + \epsilon.$$  
  Since this holds for all $\epsilon > 0$ the result
follows.   

 \end{proof}

 For consistency of notation we will denote the flat norm of a polyhedron $P$ by 
$$M_{\flat}(P) = \inf\{M(B) + M(C): P = B + \partial C, B \in \mathbf{P}_k, C \in\mathbf{P}_{k+1}\}.$$
 
\begin{prop}\label{gmt} If  
$D_j \in D_k(G) $ with $\sum E_{\flat}(D_j) < \infty$ then  a.e.
$s$ $$\sum E_{\flat}(D_j \cap H_s) <
\infty.$$
 \end{prop}
 
\begin{proof}  There exist $B_j$ and $C_j$ such that $D_j = B_j + \partial C_j$ and $$\sum E(B_j) + E(C_j) <
\infty.$$  Then
$$D_j \cap H_s = (B_j \cap H_s) - C_{js} +  \partial(C_j \cap H_s)$$ where
 $C_{js}$ is the slice   of $C_j$ as defined above.   By \ref{Qlemma}
\[
\begin{array}{rll}
E_{\flat}(D_j \cap H_s) &\le  
  E(B_j \cap H_s) + E(C_{js}) + E(C_j \cap H_s) \\&
\le
E(B_j) + E(C_{js}) + E(C_j).
\end{array}
\]
  By Lemma \ref{lintegral} $\sum\int_{-\infty}^{\infty} E(C_{js})ds \le \sum
E(C_j) < \infty.$  We conclude $\sum  E(C_{js}) <
\infty$ a.e. s.  Since $\sum  E(B_j) +  E(C_j) < \infty$ it follows that for a.e. $s$, $\sum
E_{\flat}(D_j \cap H_s) < \infty.$

\end{proof}

 Let $\Pi^0$ be a $k$-plane in $\R^3.$    If $\sigma$ is a $j$-cell in $\R^3$, $j \le k$, let $\sigma^0$ denote its
orthogonal projection into $\Pi^0.$  We may assume that $\sigma^0$ is nondegenerate, without loss of
generality.  Projection does not increase mass.  If
$S =
\sum
\delta
\sigma_i$ is a dipole
$j$-chain then $S^0 = \sum
\delta \sigma_i^0$ denotes its projection into
$\Pi^0.$  Its weight is not increased.  Furthermore, $\partial S^0 = (\partial S)^0.$   If $T = \sum \mu \tau_i$ is a
mass
$j$-chain, then $T^0 =
\sum \mu
\tau_i^0$ denotes its projection.  Again mass is not increased and $\partial T^0 = (\partial T)^0.$  If $D = S+T$ is a
dipolyhedron its orthogonal projection into $\Pi^0$ is $D^0 = S^0 + T^0$.  
 It follows that if $D = S+T$ is a $j$-dipolyhedron
then
$D^0$ is well defined $j$-dipolyhedron in $\Pi^0$ with $E(D^0) \le E(D)$ and $\partial D^0 = (\partial D)^0.$

\begin{lem} \label{WMT} Suppose $D = S+T$ is a $k$-dipolyhedron.  Let $S^0$ denote orthogonal projection of $S$
into  a
$k$-plane
$\Pi^0$ and  $T^1$  the orthogonal projection of $T$ into a $(k-1)$-plane
$\Pi^1.$ Then
$$ W(S^0)
\le E_{\flat}(D)$$ and
 $$ M(T^1) \le E_{\flat}(D)  .$$
\end{lem} 

\begin{proof}   
Let $\epsilon > 0$.  There exists a $k$-dipolyhedron $B$ and a $(k+1)$-dipolyhedron $C$ such that $D = B + \partial
C$ and
$E_{\flat}(D) > E(B) + E(C) - \epsilon.$ Write $B = S_B + T_B$ and $C = S_C + T_C$.  

For the first inequality note that since
$S_C$ has dimension
$k+1$   its projection into the $k$-plane $\Pi^0$ must be zero.    Thus $C^0 = T_C^0$.  Since $T_C$ is a
mass
$(k+1)$-chain supported in a cellular $k$-chain it follows that $(\partial T_C)^0$ is the sum of a dipole $k$-chain
$S_*$ and a mass $k$-chain
$T_*$ with $W(S_*) =  M(T_C^0) \le M(T_C) \le E(C)$.   
  
Since $  D^0 = B^0 + \partial C^0  = B^0 + S_* + T_*$ we conclude $S^0 = B^0 + S_*.$  Therefore  $$W(S^0) \le  
W(B^0) + W(S_*)
\le E(B) + E(C) < E_{\flat}(D)  + \epsilon.$$

   For the second part, note that  $(\partial T_C)^1$ vanishes since it is the
projection of a mass $(k-1)$-cycle into the $(k-1)$-plane $\Pi^1$.  Hence $D^1 = S_B^1 + T_B^1 + \partial S_C^1$
where
$\partial S_C^1$ is a dipole chain.   It follows that
$T^1 = T_B^1.$  Hence 
$$M(T^1) \le M(T_B) \le E(B) < E_{\flat}(D) + \epsilon.$$   The lemma follows since these estimates hold for every
$\epsilon > 0.$

\end{proof} 
In this section we call an open $3$-dimensional cube $Q$ {\em exceptional} relative to a sequence of
dipolyhedra $\{D_i\}$  if
$\sum E_{\flat}(D_i \cap Q) =
\infty.$  By
\ref{gmt} $Q$ is exceptional only when its faces lie on hyperplanes taken from a certain null set.  

\begin{lem}\label{S} $E_{\flat}(S) \le E_{\flat}(D).$
\end{lem} 

\begin{proof} Let
$\epsilon > 0$.  There exists $D = B + \partial C$ such that $E_{\flat}(D) > E(B) + E(C) - \epsilon.$  Then $D = S_B +
T_B +
\partial S_C + \partial T_C.$  As in \ref{WMT}, $\partial T_C = S_* + T_*$ with $W(S_*) = M(T_C).$  Thus $S = S_B +
\partial S_C + S_*.$  It follows that $$E_{\flat}(S) \le  W(S) \le W(S_B) + W(S_C) + M(T_C) \le E(B) + E(C) <
E_{\flat}(D) +
\epsilon.$$ Since this holds for each $\epsilon > 0$, we conclude  the result.
\end{proof}

\begin{thm} $E_{\flat}$ is a norm on the group $D_k(G).$
\end{thm} 

\begin{proof}  Suppose $D = S+T \ne 0 $  and  $E_{\flat}(D) = 0.$  Then either $S \ne 0$ or $T \ne 0.$   

Suppose first that
$S = 
\sum \delta
\sigma_i \ne 0 $ is nonoverlapping.   Then $\sigma_i \ne 0$ for some $i$.  By \ref{S} $E_{\flat}(S) =0.$  
Apply 
\ref{gmt} to
$D_j = S, j = 1, 2, \dots,  $    to find a nonexceptional cube $Q$, relative to $S$, such that $D \cap Q = S \cap Q = 
\delta
\sigma$ where $\sigma \ne 0$, $|\sigma| \subset |\sigma_i|$ and  $  E_{\flat}(D \cap Q) =0$.  
   By \ref{WMT} 
$$0 \ne M(\sigma) = W(\delta \sigma)   = W(S \cap Q) \le E_{\flat }(D \cap Q) = 0.$$   

Finally suppose $S = 0$ and $T = \sum \mu \tau_i \ne 0 $ is nonoverlapping.   Then $\tau_i \ne
0$ for some $i$.   By assumption
$E_{\flat}(T) = 0.$  Apply
\ref{gmt}  
  to
$D_j = T, j = 1, 2, \dots,  $    to find a  nonexceptional cube $Q$, relative to $T$,  such that  $E_{\flat}(T \cap Q) = 0
$ and 
$T \cap Q
  = 
\mu
\tau$ where $\tau \ne 0$ and
$|\tau|
\subset |\tau_i|$.   By \ref{WMT} 
$$0 \ne M(\tau) = M(\mu \tau)   = M(T \cap Q) \le E_{\flat }(D \cap Q) = 0.$$   
This contradicts the assumption that $E_{\flat}(D) = 0.$ It follows that $E_{\flat}$ is a norm.

\end{proof}
 \section{Lower semicontinuity of weight, mass and energy}
\begin{thm} \label{lweight}
Weight, mass and energy  are lower semicontinuous  in
$\mathbf{D}_k(G) $.
\end{thm}

\begin{proof}  We show that $\lim_{i \to \infty} E_{\flat}(D_i - D) = 0$ implies $$W(D) \le \liminf_{i \to \infty}
W(D_i), M(D) \le
\liminf_{i \to \infty} M(D_i), E(D) \le \liminf_{i \to \infty} E(D_i).$$ It suffices to prove these for sequences satisfying
$\sum_{i = 1}^{\infty} E_{\flat}(D_i - D) <
\infty.$    Write
$D_i = S_i + T_i$. 

Let
$Q$ be nonexceptional with respect to the sequence
$\{ D_i - D\}$ and such that $ D\cap Q = \delta \sigma$ lies in some $k$-plane $\Pi^0.$  By \ref{gmt} $\sum_{i =
1}^{\infty} E_{\flat}((D_i -D)
\cap Q ) <
\infty.$ 
   By \ref{WMT} 
$$W(((D_i - D) \cap Q)^0)= W(((S_i - S) \cap Q)^0)   \le E_{\flat}((D_i -D) \cap Q ).$$ The right
hand side tends to 0.  Therefore
$$W(D \cap Q) = W((D \cap Q)^0)  = \lim_{i \to \infty} W((D_i \cap Q)^0) \le \liminf_{i \to \infty} W(D_i \cap
Q).$$  

Given $\epsilon > 0$ there exists a finite disjoint collection of such intervals $Q$
such that $W(D \cap (Q_1 \cup \cdots \cup Q_p)^c) < \epsilon.$  From this fact and the previous inequality we deduce that $W(D) \le
\liminf_{i \to \infty} W(D_i)$.  

For the second inequality, we use \ref{gmt} to find a  nonexceptional cube $Q$ such that   $D \cap Q = S \cap Q+
\mu
\tau$ where
$\tau$ is a
$(k-1)$-cell in a
$(k-1)$-plane $\Pi^1$ and $\sum_{i = 1}^{\infty} E_{\flat}((D_i-D) \cap Q) < \infty.$  Note that $M(D \cap Q) =
M((D
\cap Q)^1 ) $ since $T\cap Q  = \mu \tau.$    By   \ref{WMT} 
$$M(((D_i - D) \cap Q)^1)= M(((T_i  -T) \cap Q)^1  ) \le  E_{\flat}((D_i -D) \cap Q).$$ The right hand side tends
to 0.  Therefore
$$M(D \cap Q) = M((D \cap Q)^1 ) = \lim_{i \to \infty} M((D_i \cap Q)^1 ) \le  
\liminf_{i
\to
\infty} M(D_i
\cap Q).$$  

Given $\epsilon > 0$ there exists a finite disjoint collection of such intervals $Q$
such that $M(D \cap (Q_1 \cup \cdots \cup Q_p)^c) < \epsilon.$  From this fact and the previous inequality we deduce
that
$M(D) \le
\liminf_{i \to \infty} M(D_i)$.  

Finally, 
\[
\begin{array} {rll}E(D) = W(D) + M(D) &\le \liminf_{i \to \infty} W(D_i) + \liminf_{i \to \infty} M(D_i) \\&\le
\liminf_{i \to \infty} (W(D_i) + M(D_i) )\\&= \liminf_{i \to \infty} E(D_i).
\end{array}
\]
\end{proof}

 We may now define the {\em energy} of a flat dipolyhedron $A $ as
 $$E(A) = \inf \{\liminf E(D_i): D_i \buildrel E_{\flat}\over\longrightarrow A \mbox{ where the } D_i \in D_k(G)
\}.$$ 
If $E(D_i) \to \infty$ for every such sequence then define $E(A)=
\infty.$   Similarly,  the {\em weight} of $A$ and {\em mass} of $A$ are defined by 
 $$W(A) = \inf \{\liminf W(D_i): D_i \buildrel E_{\flat}\over\longrightarrow A, D_i \in D_k(G) \}$$  
$$M(A) = \inf \{\liminf M(D_i): D_i \buildrel E_{\flat}\over\longrightarrow D, D_i \in D_k(G) \}.$$

 Weight    turns out   to be the area of a soap film and mass is the length of its singular curves.  Energy is their sum.

\section{Models for soap films}

\subsection*{Lipschitz dipolyhedra}   Let $U \subset \R^n$ be open.  Define $C_*(G,U)$ to be the subcomplex of
$\mathbf{D}_*(G)$ of flat dipolyhedra  satisfying $A = \lim D_i, |D_i| \subset U.$  
Suppose $f: U \to V$ is Lipschitz.    We show
that $f$ induces an operator $f_*:  C_*(G,U) \to C_*(G,V).$ 
 
\begin{lem}  If $P$ is a   polyhedron $$E_{\flat}(\delta P) \le M_{\flat}(P), \quad E_{\flat}(\mu P) \le
2M_{\flat}(P).$$
\end{lem} 
 
\begin{proof} Let $\epsilon > 0$.  There exists $P = Q + \partial R$ such that $M_{\flat}(P) > M(Q) + M(R) - \epsilon.$ 
Then
$\delta P =
\delta Q + \delta \partial R = \delta Q + \partial \delta R.$  Thus $$E_{\flat}(\delta P) \le E(\delta Q) + E(\delta R) =
M(Q) + M(R) < M_{\flat}(P)  +
\epsilon.$$    

Also, $\mu P = \mu Q + \mu \partial R = \mu Q + \partial \mu R - \delta R.$  Therefore $$E_{\flat}(\mu P)  \le E(\mu
Q) + E(\delta R) + E( \mu R) = M(Q) + 2M(R) < 2(M_{\flat}(P)  + \epsilon).$$
\end{proof}
 
 \begin{cor} \label{flatcor} If $A$ is a   flat chain then $\delta A$ and $\mu A$ are well defined flat dipolyhedra with 
$$E_{\flat}(\delta A) \le M_{\flat}(A) \mbox{ and }   E_{\flat}(\mu A) \le 2M_{\flat}(A).$$
\end{cor}

\begin{defn}[Pushforward] If $P$ is a polyhedron then
$f_*P$ is a flat chain \cite{Fleming2} and thus $f_*(\delta \sigma)  = \delta (f_*\sigma)$ and $f_*(\mu \tau) = \mu (
f_*
\tau)$ are well defined.  Extend the definition to dipolyhedra
$D = \sum \delta \sigma_i + \sum \mu \tau_i$ by linearity.   $$f_*D = f_* \left(\sum \delta \sigma_i + \sum
\mu \tau_i\right) = \sum \delta f_*\sigma_i + \sum \mu f_* \tau_i.$$  
 
\end{defn}

\begin{cor}  If $D \in D_k(G)$ then $f_*D$ is a flat dipolyhedron satisfying $\partial f_* D = f_* \partial D.$   
\end{cor} 
\begin{proof} It follows from the definitions and relations that $\partial f_* \delta \sigma = f_* \partial \delta \sigma $  
 and $\partial  f_* \mu \sigma = f_* \partial \mu
\sigma. $  Extend by linearity to dipolyhedra $D$. 
\end{proof}
If $\gamma$ is a Lipschitz Jordan curve then $ \delta \gamma$  is a
 flat dipolyhedron by \ref{flatcor}. Consider the collection $\Gamma$ of all $2$-dimensional flat
dipolyhedra spanning
$\delta \gamma$ with finite weight.    The next result shows $\Gamma$ is nonempty.  

\begin{cor}  A nondegenerate cone over a Lipschitz    Jordan curve $ {\gamma}$ supports a flat
dipolyhedron $D$ such that $\partial D =   \delta \gamma$ and $W(D) < \infty.$
\end{cor}

\begin{proof}  Assume that $\gamma$ is supported in a ball of radius $r$ and is disjoint from the origin.  Since
$\gamma$ is a flat chain it follows that
$0
\gamma$ is also a flat chain (\cite{Fleming2} p 172).    By \ref{flatcor} we know that
$\delta 0
\gamma$ is a flat dipolyhedron  with boundary  $\delta \partial 0 \gamma =   \delta \gamma$.   Furthermore,
$W(\delta 0 \gamma) =  M(0 \gamma) \le \frac{r}{2} M(\gamma).$  Set $D = \delta 0 \gamma.$
 
\end{proof}
   
 We next show that $\Gamma$ contains all
real soap films $X$ spanning $\gamma$, as observed by Plateau, letting $G = \Z_2$.  
    Consider the connected components
  $\{X_j\}$ of $X$ that are complementary to the triple branch cells $\{\tau_j\}$.   Each $X_j$ and each $\tau_j$  is
 embedded and smooth and thus flat.  Define $$D(X) = \sum \delta X_j + \sum \mu \tau_j.$$ Then $\partial D(X) =
\delta \gamma$ and
 $|D(X)| = |X|.$

\begin{ex}[Models that do not span the entire wire]  It is   easy to produce real soap films that do not span an entire
knotted wire $\gamma$.   (See Figure 3.)  Almgren called such examples striking and wrote  
\cite{A2}  \begin{quote}
 ... if one wishes to construct a theory of minimal surfaces which in particular
includes the phenomena suggested by soap films, then one must at times admit boundaries of
substantial positive thickness.  
\end{quote}   The methods of this paper support Almgren's intuition as essentially correct.  Each of our soap film
boundaries $\partial D$ is a 1-dimensional dipolyhedron bounding a two dimensional mass chain  supported in $\gamma$.   Although
$\partial D$ is algebraically closed,  there is no requirement  that the support of $\partial D$ must   itself be closed.  Figure
$3$ depicts a dipolyhedral boundary $\partial D$ that is supported in a proper arc $S$ of a Jordan curve $\gamma$.   
This dipolyhedral boundary takes the form $\partial D = \partial \mu S = \delta S - \mu \partial S;$   the mass chain $-\mu \partial S$
closes up the dipole chain $\delta S.$   Although there is a triple junction
$T$ where three surfaces of the film meet in pairs of
$120^{\circ},$ the interior of $T$ is disjoint from     $S$.   This is an example where the support of a cycle is not itself a cycle, but
an arc.

\end{ex}

\begin{figure}\label{fig.3}
\vspace{2.5in}
\hspace{-3.0in}
%\special{picture film.ps}
\begin{center} 
\resizebox{4.0in}{!}{\includegraphics*{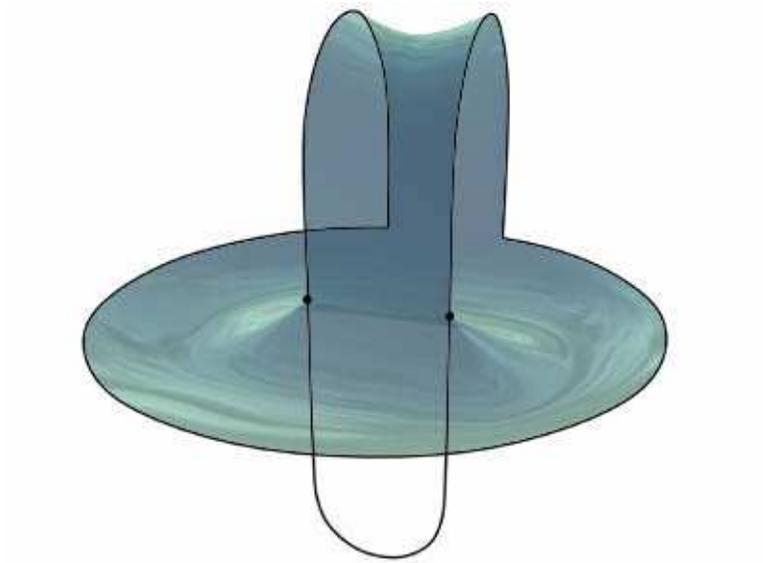}} 
\end{center}
\caption{A soap film model $D$  whose boundary is a cycle  supported in a proper arc of a Jordan curve}
\end{figure}

 \section{Future directions}
Existence and soap film regularity of an area minimizer within the set of all dipolyhedra spanning a fixed, smooth Jordan curve  remains
an open question, the solution of which would provide a  general solution to Plateau's Problem.   In a sequel
\cite{existence}, the author proves existence and regularity almost everywhere of an area minimizer, within the class of dipolyhedra with
a bound on the length of the singular curves.


\newpage \begin{thebibliography}{1}

\bibitem[A1]{A1}	Frederick J. Almgren, Jr., Existence and regularity
		almost everywhere of solutions to elliptic
		variational problems among surfaces of varying
		topological type and singularity structure,
		Ann. of Math. 87 (1968), 321-391.

\bibitem[A2]{A2} 	Frederick J. Almgren, Jr.,  Plateau's Problem,
Benjamin, New York, 1966.

\bibitem[B]{Brakke} Kenneth Brakke, Soap Films and covering spaces,
Research report, The Geometry
Center,  University of Minnesota (1993).
\bibitem[D]{Douglas}		Jesse Douglas, Solution of the
problem of Plateau,
		Trans. Amer. Math. Soc. 33 (1931), 263-321.
\bibitem[FF]{FF}		Herbert Federer and Wendell H. Fleming,
		Normal and integral currents, Ann. of Math. 72
	(1960), 458-520.


\bibitem[F]{Fleming2} Wendell H. Fleming, Flat chains over a finite
coefficient group, Trans. AMS,
Jan-Feb 1966 v 121 160-186
\bibitem[H1]{continuity} Jenny Harrison, Continuity of the integral as a function of the domain,
Journal of Geometric Analysis, 8 (1998), no. 5, 769--795. 

\bibitem[H2]{iso} Jenny Harrison, Isomorphisms of differential forms and cochains, Journal of Geometric Analysis, 8
(1998), no. 5, 797--807.


\bibitem[H3]{currents} Jenny Harrison, Geometric representations of currents and distributions, to appear, Proceedings
of Fractal Geometry and Stochastics III, Friedrichroda, Germany, 2003.

\bibitem[H4]{existence} Jenny Harrison, On Plateau's problem for soap films with a bound on energy, to appear, Journal of Geometric
Analysis.


\bibitem[H5]{discrete} Jenny Harrison, Discrete exterior calculus, preprint.


\bibitem[O]{oss} Robert Osserman, Variations on a Theme of Plateau,  Global Analysis  and
its Applications III, International Atomic Energy Agency, Vienna, 1974



\bibitem[W]{W} Hassler Whitney, Geometric Integration Theory, Princeton University Press, 1957.


 

\end{thebibliography}
\end{document}